%
\documentstyle{amsppt} 
\magnification=\magstep1
\NoRunningHeads
\NoBlackBoxes
\parindent=1em 
\vsize=7.4in


%
\topmatter
 
\title
Characterizing Hilbert space Frames with the Subframe Property            
\endtitle
\author
Peter G. Casazza 
\endauthor
\address
Department of Mathematics,
The University of Missouri,
Columbia, Missouri 65211,
USA
\endaddress
\email
pete\@casazza.math.missouri.edu\endemail
\thanks
The author was supported by NSF DMS-9201357, Danish NSRC grant 9401958, Missouri Research Board grant C-3-41743, and a Missouri Research Council Summer Fellowship.
\endthanks
%
\abstract
We characterize frames which satisfy the subframe property (i.e. frames with the property that every subset is a frame for its closed linear span).  With this characterization we can answer most of the problems from the literature concerning the subframe property, including its relationship to Riesz frames and the projection methods.
\endabstract
\endtopmatter
\document
\baselineskip=15pt

\heading{1. Introduction}
\endheading
\vskip10pt

A sequence $(f_{i})_{i=1}^{\infty}$ in a Hilbert space $H$ which is a frame for its closed linear span is called a {\bf frame sequence}.  If every subsequence of $(f_{i})_{i=1}^{\infty}$ is a frame sequence, we say that the frame has the {\bf subframe property}. If $(f_{i})_{i=1}^{\infty}$ is a frame for $H$ with the subframe property and additionally there are uniform upper and lower frame bounds for all subsequences of the frame, then we call $(f_{i})_{i=1}^{\infty}$ a {\bf Riesz frame}. 
Riesz frames were introduced in \cite{6} where it was shown that every Riesz frame for $H$ contains a subset which is a Riesz basis for $H$.  The projection methods \cite{4} play a central role in evaluating truncation error which arises in computing approximate solutions to moment problems, as well as handling the very difficult problem of computing dual frames.  There were many natural questions arising from the literature concerning the interrelationships between Riesz frames, frames with the subframe property, and the projection methods \cite{2,4,5,6,8}.  In this paper we characterize Riesz frames and frames with the subframe property which allows us to answer most of these questions.

\heading{2.  Riesz Frames}
\endheading
\vskip10pt

If ${\Cal F}$ is a subset of $H$, we write {\bf span ${\Cal F}$} for the closed linear span of ${\Cal F}$.
A sequence $(f_{i})_{i=1}^{\infty}$ in $H$ is called a {\bf frame} for $H$ if there are positive constants $A,B$ satisfying,
$$
A\|f\|^{2} \le \sum_{i=1}^{\infty}|<f,f_{i}>|^{2} \le B\|f\|^{2}, \ \ \ \ \forall f\in H.  \tag 2.1
$$
We call $A,B$ the lower and upper frame bounds respectively.  In general, a subset of a frame need not be a frame for its closed linear span.  But clearly $B$ is an upper frame bound for every subset of the frame (i.e.  It is only the lower frame bound that might be lost when switching to a subset of a frame).  For a Riesz frame, the common frame bounds for all subsets of the frame will be called the {\bf Riesz frame bounds}.  The largest $A$ and the smallest $B$ satisfying (2.1) are called the {\bf optimal frame bounds}.    An unconditional basis $(f_{i})_{i\in I}$ for $H$ is called a {\bf Riesz basis}.  Equivalently, $(f_{i})_{i\in I}$ is a Riesz basis if it is total and there are constants $c,C$ so that for every sequence of scalars $(a_{i})_{i\in I}$ we have
$$
c\sqrt{\sum_{i\in I}|a_{i}|^{2}} \le \|\sum_{i\in I}a_{i}f_{i}\| \le C\sqrt{\sum_{i\in I}|a_{i}|^{2}}. \tag 2.2
$$
The largest $c$ and the smallest $C$ satisfying (2.2) are called the {\bf Riesz basis constants} for $(f_{i})_{i\in I}$.  If $(f_{i})_{i\in I}$ is a Riesz basis, then \cite{7} the Riesz basis constants equal the square root of the optimal frame bounds.  Finally, we say that two frames $(f_{i})_{i=1}^{\infty}, (g_{i})_{i=1}^{\infty}$ are {\bf equivalent} if there is an isomorphism $T:H\rightarrow H$ with $T(f_{i}) = g_{i}$, for all $i = 1,2,\ldots$.

We start with an elementary observation concerning Riesz frames. 

\proclaim{Proposition 2.1}
For a frame $(f_{i})_{i=1}^{\infty}$ for $H$, the following are equivalent:

(1)  $(f_{i})_{i=1}^{\infty}$ is a Riesz frame, 

(2)  There is an $A > 0$ so that for every finite set of natural numbers $\Delta$ for which $(f_{i})_{i\in \Delta}$ is linearly independent, the family $(f_{i})_{i\in \Delta}$ has lower Riesz basis bound $A$.
\endproclaim

\demo{Proof}
$\Rightarrow$  If $(f_{i})_{i\in \Delta}$ is linearly independent, then the lower Riesz basis constant for this set equals the square root of the lower Riesz frame bound.

$\Leftarrow$  It is only the lower frame bound that needs to be checked.  For any finite set of natural numbers $\Gamma$, let $(f_{i})_{i\in \Delta}$ be a maximal linearly independent subset, where $\Delta \subset \Gamma$.  Then the lower frame bound of $(f_{i})_{i\in \Gamma}$ is greater than or equal to the lower frame bound of $(f_{i})_{i\in \Delta}$ which is equal to the square root of the lower Riesz basis constant, $\sqrt{A}$.  So $(f_{i})_{i=1}^{\infty}$ is a Riesz frame. 
\enddemo

This remark yields a short proof of a result of Christensen \cite{6}.

\proclaim{Corollary 2.2 (Christensen)}
Every Riesz frame contains a Riesz basis.
\endproclaim

\demo{Proof}
Choose a maximal linearly independent subset of the frame.  This is a Riesz basis, by Proposition 2.1.
\enddemo

We now introduce some of the notation which will be used throughout the paper.  If $(g_{i})_{i\in I}$ is a Riesz basis for $H$, and $\Delta \subset I$, we let $P_{\Delta}$ be the natural projection of $\text{span}(g_{i})_{i\in I}$ onto $\text{span}(g_{i})_{i\in \Delta}$.  That is, $P_{\Delta}\sum_{i\in I}a_{i}g_{i} = \sum_{i\in \Delta}a_{i}g_{i}$.  We will also write $P_{n} = P_{\{1,2,\ldots,n\}}$, and for $m < n$, $P_{n,m} = P_{n} - P_{m}$.  If $(f_{i})_{i\in I}$ is a frame with frame bounds $A,B$, and $P$ is an orthogonal projection on $H$, then $(Pf_{i})_{i\in I}$ is a frame sequence with frame bounds $A,B$.  Conversely, if $(f_{i})_{i\in I}$ (respectively, $(g_{j})_{j\in \Gamma}$) is a frame for $P(H)$ (respectively $(I - P)(H)$) with frame bounds $A_{1},B_{1}$ (respectively, $A_{2}, B_{2}$), then $((f_{i})_{i\in I},(g_{j})_{j \in \Gamma})$ is a frame for $H$ with frame bounds,
$$
A = \text{min}\{A_{1},A_{2}\}, \ \ \ \  B = \text{max}\{B_{1},B_{2}\}.
$$
We will make extensive use of a slight extension of these properties which we now state.

\proclaim{Proposition 2.3}
Let $(f_{i})_{i=1}^{\infty}$ be a sequence in $H$ with upper frame bound $B$.  Let $\Delta$ be a subset of the natural numbers and $P$ denote the orthogonal projection of $H$ onto $\text{span}(f_{i})_{i\in \Delta}$.

(1)  If $(f_{i})_{i\in \Delta}$ is a frame with frame bounds $A_{1},B$, and $((I - P)f_{i})_{i \in {\Delta}^{c}}$ is a frame sequence with frame bounds $A_{2},B$, then $(f_{i})_{i=1}^{\infty}$ is a frame for $H$ with frame bounds $\frac{A_{1}A_{2}}{8B},B$.

(2)  If $(f_{i})_{i=1}^{\infty}$ is a frame with frame bounds $A,B$  then $((I - P)f_{i})_{i\in {\Delta}^{c}}$ is a frame sequence with frame bounds $A,B$.   
\endproclaim

\demo{Proof}
(1)  For any $f\in H$ we have,
$$
\sum_{i=1}^{\infty}|<f,f_{i}>|^{2} = \sum_{i\in \Delta}|<f,f_{i}>|^{2} + \sum_{i\in {\Delta}^{c}}|<f,f_{i}>|^{2}  \tag 2.3
$$
$$
= \sum_{i\in \Delta}|<Pf,f_{i}>|^{2} + \sum_{i\in {\Delta}^{c}}|<Pf,Pf_{i}> + <(I - P)f,(I - P)f_{i}>|^{2} 
$$
$$
\ge A_{1}\|Pf\|^{2} + \left[\sqrt{\sum_{i\in {\Delta}^{c}}|<(I - P)f,(I - P)f_{i}>|^{2}} - \sqrt{\sum_{i\in {\Delta}^{c}}|<Pf,Pf_{i}>|^{2}}\right]^{2}
$$
$$
\ge A_{1}\|Pf\|^{2} + \left[\sqrt{A_{2}}\|(I - P)f\| - \sqrt{B}\|Pf\| \right]^{2} 
$$
Now, there are two possibilities.

\proclaim{Case I}
$\|Pf\|^{2} \ge \frac{A_{2}}{8B}\|f\|^{2}$ 
\endproclaim
In this case, inequality (2.3) and the fact that $\frac{A_{2}}{B} \le 1$ immediately yields,
$$
\sum_{i=1}^{\infty}|<f,f_{i}>|^{2} \ge  \frac{A_{1}A_{2}}{8B}\|f\|^{2}.
$$

\proclaim{Case II}
$\|Pf\|^{2} \le  \frac{A_{2}}{8B}\|f\|^{2}$ 
\endproclaim
In this case, since $\frac{A_{2}}{8B} \le \frac{1}{8}$, we have that $\|(I - P)f\|^{2} \ge \frac{1}{2}$. This combined with inequality (2.3) yields,
$$
\sum_{i=1}^{\infty}|<f,f_{i}>|^{2} \ge \left[\sqrt{A_{2}}\|(I - P)f\| - \sqrt{B}\|Pf\| \right]^{2} 
$$
$$
\ge \left[\sqrt{\frac{A_{2}}{2}}\|f\| - \sqrt{\frac{A_{2}}{4}}\|f\|\right]^{2} = \frac{A_{2}}{8}\|f\|^{2} \ge \frac{A_{1}A_{2}}{8B}\|f\|^{2}.
$$

(2)  By our assumptions, $(I - P)H = \text{span}((I - P)f_{i})_{i\in {\Delta}^{c}}$.  So for any $f \in (I - P)H$,
$$
A\|f\|^{2} \le \sum_{i=1}^{\infty}|<f,f_{i}>|^{2} = \sum_{_i \in {\Delta}^{c}}|<f,(I - P)f_{i}>|^{2}.
$$
\enddemo

We next give a (slightly internal) classification of Riesz frames which is of some interest itself, and will be important later for our classification of frames with the subframe property.

\proclaim{Theorem 2.4}
The following are equivalent for a frame $(f_{i})_{i=1}^{\infty}$:

(1)  $(f_{i})_{i=1}^{\infty}$ is a Riesz frame,

(2)  $(f_{i})_{i=1}^{\infty}$  can be divided into two subsets, $(g_{i})_{i=1}^{\infty},(h_{i})_{i\in \Gamma}$ satisfying,

\ \ \ \ \ (i)  $(g_{i})_{i=1}^{\infty}$ is a Riesz basis for $H$,

\ \ \ \ (ii)  There is an $A_{0} > 0$ so that for each subset $\Delta$ of the natural numbers, and ${\Gamma}_{1}\subset \Gamma$,

\ \ \ \  the set $(P_{\Delta}h_{i})_{i\in {\Gamma}_{1}}$ is a frame sequence with lower frame bound $A_{0}$.

Moreover, in this case, if $A,B$ are the Riesz frame bounds for $(f_{i})_{i=1}^{\infty}$, then there is natural number k so that we can write 
$h_{i} = \sum_{j\in {\Delta}_{i}}h_{i}(j)g_{j}$, with $|{\Delta}_{i}| \le k$, and 
$A^{2} \le |h_{i}(j)| \le B^{2}, \forall j\in {\Delta}_{i}$.
\endproclaim

\demo{Proof}
$(1) \Rightarrow (2)$  Since $(f_{i})_{i=1}^{\infty}$ is a Riesz frame, by Corollary 2.2, it contains a Riesz basis, say $(g_{i})_{i=1}^{\infty}$.  Let $(h_{i})_{i\in \Gamma}$ be the remaining elements of the frame, and assume that $A,B$ are the Riesz frame bounds for $(f_{i})_{i=1}^{\infty}$.  It suffices to prove the theorem for any frame equivalent to our frame.  So, by taking the natural isomorphism of $(g_{i})_{i=1}^{\infty}$ to an orthonormal basis for $H$, we may assume, without loss of generality, that $(g_{i})_{i=1}^{\infty}$ is an orthonormal basis for $H$.  However, the Riesz frame bounds have to be adjusted by the norm of the isomorphism to $A^{2}, B^{2}$.  If $\Delta,{\Gamma}_{1}$ be as in (2)(ii), then $((g_{i})_{i\in {\Delta}^{c}},(h_{i})_{i\in {\Gamma}_{1}})$ is a frame sequence with frame bounds $A^{2},B^{2}$, and $(g_{i})_{i=1}^{\infty}$ is an orthonormal basis.  If $P_{{\Delta}^{c}}$ is the natural projection of $H$ onto $(g_{i})_{i\in {\Delta}^{c}}$, then $I - P_{{\Delta}^{c}} = P_{\Delta}$.  Proposition 2.3(2) now yields that $(P_{\Delta}h_{i})_{i\in {\Gamma}_{1}}$ is a frame sequence with frame bounds $A^{2},B^{2}$.  This concludes the proof that $(1)$ implies $(2)$.  To check the "moreover" part, write $h_{i} = \sum_{j\in {\Delta}_{i}}h_{i}(j)g_{j}$, where $h_{i}(j) \not= 0$, for $j\in {\Delta}_{i}$.  For any $i = 1,2,\ldots$,and any $j\in {\Delta}_{i}$, consider the subset $F = \{h_{i}\} \cup \{g_{m}: j\not= m \in {\Delta}_{i}\}$.  Then $g_{j}\in \text{span}F$ and this set has frame bounds $A^{2},B^{2}$ implies,
$$
\sum_{j\not= m \in {\Delta}_{i}}|<g_{m},g_{j}>|^{2} + |<h_{i},g_{j}>|^{2} = |h_{i}(j)|^{2} \ge A^{2}.  \tag 2.4 
$$ 
Also, 
$$
|h_{i}(j)|^{2} \le \|h_{i}\|^{2} \le B^{2}.  \tag 2.5
$$      
Since $\text{sup}_{1\le i < \infty}\|f_{i}\| < \infty$, the existence of $k$ is now immediate from (2.4) and (2.5).

$(2) \Rightarrow (1)$  Let $(f_{i})_{i=1}^{\infty} = ((g_{i})_{i=1}^{\infty},(h_{i})_{i \in \Gamma})$ be a sequence of vectors in $H$ satisfying (2).  Again we can start by taking the natural isomorphism of $(g_{i})_{i=1}^{\infty}$ onto an orthonormal basis $(e_{i})_{i=1}^{\infty}$.  This will change the $A_{0}$ in (2)(ii) to say $A$.  Letting $\Delta$ equal the natural numbers and ${\Gamma}_{1} = \Gamma$ in (2)(ii), we see that $(h_{i})_{i\in \Gamma}$ is a frame with frame bounds $A,B$.  So $(f_{i})_{i=1}^{\infty}$ has a finite upper frame bound $1 + B$.  Choose a subset of our set of vectors of the form: $((g_{i})_{i \in {\Delta}}, (h_{i})_{i\in {\Gamma}_{2}})$.  Let ${\Gamma}_{1} = \{i\in {\Gamma}_{2}: P_{\Delta}h_{i} \not= 0\}$.  By our assumption (2)(ii), $(P_{{\Delta}^{c}}h_{i})_{i\in {\Gamma}_{1}}$ has lower frame bound $A$.  Applying Proposition 2.3 (2) (recall that $(g_{i})_{i=1}^{\infty}$ is an orthonormal basis) we have that $(f_{i})_{i=1}^{\infty}$ has lower frame bound $\frac{A}{8}$.  So $(f_{i})_{i=1}^{\infty}$ is a Riesz frame.
\enddemo

Let us recall some notation.  If $(f_{i})_{i=1}^{\infty}$ is a basis for its span, we say that a sequence $(g_{i})_{i=1}^{\infty}$ is {\bf disjointly supported} with respect to $(f_{i})_{i=1}^{\infty}$ if there exist a disjoint family of subsets of the natural numbers $({\Delta}_{i})_{i=1}^{\infty}$ so that
$$
g_{i} \in \text{span}(f_{j})_{j\in {\Delta}_{i}}, \ \ \ \ \forall i.
$$
That is, the supports of the $g_{i}$, relative to the basis $(f_{i})_{i=1}^{\infty}$, are disjoint.

Theorem 2.4 shows that Riesz frames have a somewhat exact form.  The next corollary gives a further restriction on Riesz frames.  
  
\proclaim{Corollary 2.5}
Every Riesz frame for $H$ is equivalent to one of the form 
$$
((e_{i})_{i=1}^{\infty},(f_{ij})_{i=1,j=1}^{\ k \ ,\ {\infty}})
$$ 
where $(e_{i})_{i=1}^{\infty}$ is a orthonormal basis for $H$, for each $1\le i \le k$, $(f_{i,j})_{j=1}^{\infty}$ is disjointly supported with respect to $(e_{i})_{i=1}^{\infty}$, and the non-zero coordinates (with respect to the orthonormal basis $(e_{i})$)  satisfy $A \le |f_{i,j}(n)| \le B$ for some $A,B > 0$, and there is a natural number $K$ so that
$$
|\{n : f_{i,j}(n) \not= 0\}| \le K.
$$ 
\endproclaim

\demo{Proof}
Let $B_{0} $ be the upper Riesz frame bound for $(f_{i})_{i=1}^{\infty}$, and choose a natural number $K$ so that $\frac{K}{k}(B)^{2} > B_{0}$.
Basically, we will apply the pigeonhole principle to $(h_{i})$ in Theorem 4.4 to divide it into at most K-sets, $G_{1},G_{2},\ldots,G_{K}$ where the $h_{i}$ in $G_{j}$ are disjointly supported.  We start by putting $h_{1}$ into $G_{1}$.  If $h_{2}$ has disjoint support from $h_{1}$, put it also into $G_{1}$, otherwise, put it in $G_{2}$.  We continue by induction.  Assume that $h_{1},h_{2},\ldots,h_{n}$ have been distributed into the sets so that the elements of each set are disjointly supported.  If $h_{n+1}$ is disjoint from all the elements of $G_{1}$, put it in $G_{1}$.  If not, go to $G_{2}$ and so on.  If we reach set $G_{K}$, then by assumption, $h_{n+1}$ has a non-zero coordinate in common with at least one element from each of the sets $G_{1},G_{2},\ldots,G_{K-1}$.  But, by Theorem 4.4, $h_{n+1}$ has only k non-zero coordinates.  Hence, $h_{n+1}$ has a fixed non-zero coordinate, say m, in common with $\frac{K}{k}$ of the $h_{i}$.  Hence,
$$
\sum_{i}|<e_{m},h_{i}>|^{2} \ge \frac{K}{k}(B)^{2} > B_{0},
$$
which is a contradiction.  Thus, $h_{n+1}$ must go into at least one of the sets.   
\enddemo

The next corollary shows that Corollary 2.5 comes close to classifying Riesz frames (All we are missing in Corollary 2.5 is condition (2) of Corollary 2.6).

\proclaim{Corollary 2.6}
Let $A,B > 0$, and $K$ be a natural number.  Let $(e_{i})_{i=1}^{\infty}$ be an orthonormal basis for $H$, and $(f_{ij})_{i=1,j=1}^{\ k \ \ \ {\infty}}$ be vectors in $H$ satisfying: 

(1)  The non zero coordinates of $f_{i,j}$ (with respect to the orthonormal basis $(e_{i})$)  satisfy:

\ \ \ \ \ (i)   $A \le |f_{i,j}(n)|^{2} \le B$,

\ \ \ \ (ii)  $|\{n : f_{i,j}(n) \not= 0\}| \le K.$

(2)  $\text{span}(f_{ij})_{j=1}^{\infty} \subset \text{span}(f_{i-1,j})_{j=1}^{\infty}$, $\forall 2\le i \le k$,

(3)  Each $(f_{ij})_{j=1}^{\infty}$ is a disjointly supported sequence with respect to $(f_{i-1,j})_{j=1}^{\infty}$ (with

\ \ \ \ \ $f_{0,j} = g_{j}$, for all $j = 1,2,\ldots$).

Then $((e_{i})_{i=1}^{\infty},(f_{ij})_{i=1,j=1}^{\ k \ ,\ \ {\infty}})$ is a Riesz frame for $H$ with Riesz frame bounds
$$
\frac{1}{D^{k}8^{k}\prod_{i=1}^{k}(1 + iD)},\ \ \ 1 + kD,
$$
where $D = \frac{KB}{A}$.           
\endproclaim

\demo{Proof}
We will do the proof in three steps. 

\proclaim{Step I}
We start with a calculation.
\endproclaim

Let $\Delta$ be a subset of the natural numbers and $P_{\Delta}$ denote the natural projection of $H$ onto $\text{span}(e_{i})_{i\in \Delta}$. By deleting the $f_{i,j}$ with support in $\Delta$ and reindexing, we may assume that $P_{{\Delta}^{c}}f_{i,j} \not= 0$, for all $1 \le i \le k$, and $j=1,2,\ldots$.
We will work with the family $((e_{i})_{i\in \Delta},(f_{i,j})_{i=1,j=1}^{\ k \ \ \ \infty})$.  Fix $2\le i_{0}\le k$ and note that $(P_{{\Delta}^{c}}f_{i_{0},j})_{j=1}^{\infty}$ is an orthogonal sequence in $H$ with $A \le \|f_{i_{0},j}\|^{2} \le KB$.  By taking the natural isomorphism
$$
T(P_{{\Delta}^{c}}f_{i_{0},j}) = \frac{P_{{\Delta}^{c}}f_{i_{0},j}}{\|P_{{\Delta}^{c}}f_{i_{0},j}\|},
$$
we have that $\sqrt{A} \le \|T\| \le \sqrt{KB}$.  For $i_{0} \le i \le k$, let $g_{i,j} = T(P_{{\Delta}^{c}}f_{i,j})$.  It follows that, 
\vskip10pt
(2.6)  $(g_{i_{0},j})_{j=1}^{\infty}$ is an orthonormal basis for its closed linear span 
\vskip10pt

(2.7)  $\text{span}(g_{i,j})_{j=1}^{\infty} \subset \text{span}(g_{i-1,j})_{j=1}^{\infty}$, 
\vskip10pt
(2.8)  Each $(g_{i,j})_{j=1}^{\infty}$ is a disjointly supported sequence with respect to $(g_{i-1,j})_{i=1}^{\infty}$.
\vskip10pt
kFor all $i_{0}\le i \le k$, we can choose subsets of the natural numbers ${\Delta}_{i,j}$ so that
$$
g_{i,j} = \sum_{m\in {\Delta}_{i_{0},j}}a_{m}f_{1,m}, \ \ \ a_{m}\not= 0, \forall m\in {\Delta}_{i,j}.  \tag 2.9
$$
By our assumption (1)(i), the non-zero coordinates of $P_{{\Delta}^{c}}f_{i,j}$ (relative to the Riesz basis $(e_{i})_{i=1}^{\infty}$) satisfy:
$$
A \le |P_{{\Delta}^{c}}f_{i,j}(n)|^{2} \le B.  \tag 2.10
$$
Therefore,
$$
\frac{A}{KB} \le |g_{i,j}(n)|^{2} \le \frac{KB}{A}.  \tag 2.11
$$
By (2.9) and (2.11) we have,
$$
\frac{A}{KB} \le |a_{m}f_{1,m}(n)|^{2} \le \frac{KB}{A}.  \tag 2.12
$$
Combining (2.11) and (2.12) we have
$$
\frac{1}{D} = \frac{A}{KB} \le |a_{m}| \le \frac{KB}{A} = D.  \tag 2.13
$$
It follows that for $2\le i_{0} \le i \le k$, and for every $j$, the non-zero coordinates $g_{i,j}(m)$ (This denotes the coordinates of $g_{i,j}$ with respect to the orthonormal basis $(g_{i_{0},j})_{j=1}^{\infty}$) satisfy
$$
\frac{1}{D} \le |g_{i,j}(m)| \le D, \tag 2.14
$$
Also, note that the number of these non-zero coordinates is still $\le K$.

We will prove the corollary by induction on $k$ with the hypotheses of the corollary except that we will assume that our family satisfies (2.14) and replacing $A,B$ in assumption (1)(i) by $\frac{1}{D}, D$ respectively.

\proclaim{Step II}
Starting the induction.  i.e.  The case $k=1$.
\endproclaim

Since $(e_{i})_{i=1}^{\infty}$ is an orthonormal basis for $H$, the $(f_{1,j})_{j=1}^{\infty}$ are disjointly supported, and $\frac{1}{D} \le \|f_{1,j}\| \le D$, it follows that $(f_{1,j})_{j=1}^{\infty}$ has Riesz basis constants $\sqrt{\frac{1}{D}},\sqrt{D}$, and hence frame bounds $\frac{1}{D},D$.  So $((e_{i})_{i=1}^{\infty},(f_{1,j})_{j=1}^{\infty})$ has upper frame bound $\le 1 + D$.  Let $((e_{i})_{i\in \Delta},(f_{1,j})_{j\in \Gamma})$ be a subset of our set of vectors.   Let $P_{\Delta}$ be the natural projection of $H$ onto $\text{span}(g_{i})_{i\in \Delta}$. and let
$$
\Lambda = \{j\in \Gamma:P_{{\Delta}^{c}}f_{1,j}\not= 0\}.
$$
Now, $P_{\Delta}f_{1,j} = f_{1,j}$, for all $j\in \Gamma - \Lambda$.  So $((e_{i})_{i\in \Delta},(f_{1,j})_{j\in \Gamma - \Lambda})$ is a frame with frame bounds $1, 1 + D$.   Now, $(P_{{\Delta}^{c}}f_{1,j})_{j\in \Lambda}$ is a disjointly supported sequence of vectors with respect to $(e_{i})_{i \in {\Delta}^{c}}$ for which: $\frac{1}{D} \le \| P_{{\Delta}^{c}}f_{1,j}\|^{2} \le D$.  Hence, this is a Riesz basis with constants $\sqrt{\frac{1}{D}},\sqrt{D}$ and lower frame bound $A\frac{1}{D}$.  By Proposition 2.3 (1), it follows that
$((e_{i})_{i\in \Delta},(f_{1,j})_{j\in \Gamma})$ is a frame with frame bounds $\frac{1}{D8(1 + D)}, 1 + D$.  So our family is a Riesz frame with the bounds specified in the corollary.  

\proclaim{Step III}
The induction step.
\endproclaim

Assume the result holds for some $k-1$, and we will prove that it holds for $k$.   Choose a subfamily of our set given by:  $((e_{i})_{i\in \Delta},(f_{i,j})_{i=1,j\in {\Delta}_{i}}^{\ k})$.  For each $1\le i \le k$, $(f_{i,j})_{j=1}^{\infty}$ is an orthogonal sequence satisfying $\frac{1}{D} \le \|f_{i,j}\|^{2} \le D$, and so this family has upper frame bound D.  Hence, $((e_{i})_{i\in \Delta},(f_{i,j})_{i=1,j=1}^{\ k\ \ \infty})$ has upper frame bound $1 + kD$. Since $(e_{i})_{i\in \Delta}$ is an orthonormal sequence, $((e_{i})_{i\in \Delta},(P_{\Delta}f_{i,j})_{i=1,j=1}^{\ k\ \ \infty}$ has frame bounds $1, 1 + kD$.  So, without loss of generality, we may assume that $P_{{\Delta}^{c}}f_{i,j}\not= 0$, for all $j\in {\Delta}_{i}$.  Let $i_{0} = 1$ in Step I to obtain the corresponding $g_{i,j}$.     
By Step I, we can apply the induction hypothesis to the family $((g_{1,j})_{j=1}^{\infty},(g_{i,j})_{i=2,j=1}^{\ \ k\ \ \infty})$ to discover that this is a Riesz frame with Riesz frame bounds,
$$
\frac{1}{D^{k-1}8^{k-1}\prod_{i=1}^{k-1}(1 + iD)},\ \ \ 1 + (k-1)D. \tag 2.15
$$

That is, $(TP_{{\Delta}^{c}}f_{i,j})_{i=1,j\in {\Delta}_{i}}^{\ k})$ is a frame with frame bounds given by (2.15).  Therefore, $(P_{{\Delta}^{c}}f_{i,j})_{i=1,j\in {\Delta}_{i}}^{\ k})$ is a frame with lower frame bound
$$
\frac{1}{D^{k}8^{k-1}\prod_{i=1}^{k-1}(1 + iD)}.
$$

Applying Proposition 2.3 (1), we see that $((e_{i})_{i\in \Delta},(f_{i,j})_{i=1,j\in {\Delta}_{i}}^{\ k})$ is a frame with frame bounds
$$
\frac{1}{D^{k}8^{k}\prod_{i=1}^{k}(1 + iD)},\ \ \ 1 + kD.
$$
  This proves that our original family is a Riesz frame with the stated frame bounds, and concludes the proof of corollary 2.6       

\enddemo

\heading{3.  Characterizing Frames with the Subframe Property}
\endheading
\vskip10pt

In this section we characterize of frames having the subframe property. 
To simplify the proof of the theorem, we first make an elementary observation.

\proclaim{Lemma 3.1}
If $(f_{i})_{i=1}^{\infty}$ is a frame for $H$, $G$ is a finite dimensional subspace of $H$ and $P$ is the orthogonal projection of $H$ onto $G$, then
$$
\sum_{i=1}^{\infty}\|Pf_{i}\|^{2} < \infty.
$$
\endproclaim

\demo{Proof}
Let $\{e_{1},e_{2},\ldots,e_{n}\}$ be an orthonormal basis for $G$.  Then
$$
\sum_{i=1}^{\infty}\|Pf_{i}\|^{2} = \sum_{i=1}^{\infty}\sum_{j=1}^{n}|<Pf_{i},e_{j}>|^{2} =
\sum_{i=1}^{\infty}\sum_{j=1}^{n}|<f_{i},Pe_{j}>|^{2}
$$
$$
= \sum_{i=1}^{\infty}\sum_{j=1}^{n}|<f_{i},e_{j}>|^{2}
= \sum_{j=1}^{n}\sum_{i=1}^{\infty}|<f_{i},e_{j}>|^{2}
\le \sum_{j=1}^{n}B\|e_{j}\|^{2} = nB.
$$
\enddemo

Now we are ready to prove the main theorem of this paper.

\proclaim{Theorem 3.2}
For a frame $(f_{i})_{i=1}^{\infty}$ the following are equivalent:

(1) $(f_{i})_{i=1}^{\infty}$ has the subframe property,

(2)  The frame $(f_{i})_{i=1}^{\infty}$ can be divided into three sets of vectors, $(g_{i})_{i=1}^{\infty},(h_{i})_{i\in \Gamma}, (k_{i})_{i=1}^{n}$ where $\Gamma$ may be finite or infinite, $(g_{i})_{i=1}^{\infty}$ is a Riesz basis for $H$, and there is a natural number $m$ so that if $G = \text{span}(g_{i})_{i=1}^{m}$, then $h_{i}$ is of the form $h_{i} = h_{i}^{1} + h_{i}^{2}$ with $h_{i}^{2}\in G, h_{i}^{1}\in G^{\perp}$ and satisfying:

\ \ \ \ \ (i)  The $k_{i}$ have infinite support, 

\ \ \ \ (ii)  $\sum_{i\in \Gamma}\|h_{i}^{2}\|^{2} < \infty$,

\ \ \ (iii)  $((g_{i})_{i=1}^{\infty},(h_{i}^{1})_{i\in \Gamma})$ is a Riesz frame for $H$.    
\endproclaim

\demo{Proof}
$(1)\Rightarrow (2)$  By Casazza, Christensen \cite{3}, $(f_{i})_{i=1}^{\infty}$ contains a Riesz basis, say $(g_{i})_{i=1}^{\infty}$.  To simplify the proof, we take the natural isomorphism of $(g_{i})_{i=1}^{\infty}$ to an orthonormal basis $(e_{i})_{i=1}^{\infty}$ and see that, without loss of generality, we may assume that $(g_{i})_{i=1}^{\infty}$ is an orthonormal basis for $H$.  Let $(k_{i})_{i\in \Lambda}$ be the elements of $(f_{i})_{i=1}^{\infty}$ with infinite support with respect to $(g_{i})_{i=1}^{\infty}$, and let $(h_{i})_{i\in \Gamma}$ be the remaining elements of the frame (i.e.  The elements of the frame which are not one of the $g_{i}$ and which have finite support).  We can write:
$$
h_{i} = \sum_{j\in {\Omega}_{i}}h_{i}(j)g_{j},
$$
where $|{\Omega}_{i}| < \infty$, and $h_{i}(j) \not= 0, \forall j\in {\Omega}_{i}$.   

\proclaim{Step I}
$|\Lambda| < \infty$. 
\endproclaim

We proceed by way of contradiction.  So assume we have infinitely many infinitely supported vectors $(k_{i})_{i=1}^{\infty}$.  We must construct a subset of our frame which is not a frame for its closed linear span.  To do this, we apply an inductive construction to the two conditions below:
$$
\|k_{i}\|^{2} = \sum_{j=1}^{\infty}|k_{i}(j)|^{2} < \infty, \ \ \ \ \forall i=1,2,\cdots,
$$
$$
\sum_{i=1}^{\infty}|<k_{i},g_{j}>|^{2} = \sum_{i=1}^{\infty}|k_{i}(j)|^{2} \le B, \ \ \ \ \forall j=1,2,\cdots.
$$
By alternately applying these two conditions and induction, we can find sequences of natural numbers $i_{1} < i_{2} < i_{3} < \cdots$ and $j_{1} < j_{2} < j_{3} < \cdots$ so that
$$
0 < \sum_{n=1}^{\infty}|k_{i_{n}}(j_{m})|^{2} < \frac{1}{m}, \ \ \ \ \forall m = 1,2,3,\ldots. \tag 3.1
$$
We will sketch the beginning of this induction proof.  From the first condition, we can choose a $i_{1} = 1, j_{1}$ so that:

$$
0 < |k_{i_{1}}(j_{1}|^{2} < \frac{1}{2}.
$$
The second condition allows us to switch to a subsequence of $(k_{i})_{i=1}^{\infty}$, starting with $k_{i_{1}}$ (call it $(k_{i})_{i=i_{1}}^{\infty}$) so that
$$
\sum_{n=i_{1}+1}^{\infty}|k_{n}(j)|^{2} < \frac{1}{2}.
$$
Now, using the first condition, we can find a natural number  m so that $|k_{i_{1}}(j)|^{2} < \frac{1}{(3)2}$, for all $j\ge m$.  Choose any $i_{2} > i_{1}$.  Since $k_{i_{2}}$ is infinitely supported, there is some $j_{2} > j_{1}$ so that
$$
0 \not= |k_{i_{2}}(j_{2})|^{2} < \frac{1}{(3)(2)}.
$$
By condition 2 again, we can choose an $i_{2} > i_{1}$ and a subset of the $(k_{i})_{i=i_{1}}^{\infty}$ (denote it $(k_{i_{1}},k_{i_{2}},k_{i_{2}+1},k_{i_{2}+2},\ldots)$)) so that
$$
\sum_{n=i_{2}+1}^{\infty}|k_{n}(j_{2}|^{2} \le \frac{1}{(3)(2)}.
$$
Now choose any $i_{3} > i_{2}$ and a natural number m so that
$$
\sum_{n=1}^{2}|k_{i_{n}}(j)|^{2} < \frac{1}{(3)(3)}, \ \ \ \ \forall j\ge m.
$$
Again, since $k_{i_{3}}$ is infinitely supported, there is some $j_{3} > j_{2}$ so that,
$$
0 < |k_{i_{3}}(j_{2})|^{2} < \frac{1}{(3)(3)}.
$$
and by switching to a subsequence of $(k_{i})$ we may assume that
$$
\sum_{n=i_{3}+1}^{\infty}|k_{n}(j_{2})|^{2} < \frac{1}{(3)(3)}.
$$
Now continue this construction by induction.  

Finally, let $\Delta = \{j_{m}:m=1,2,3,\ldots\}^{c}$ and consider the subframe of our frame given by: $((g_{i})_{i\in \Delta},(k_{i_{n}})_{n=1}^{\infty})$.  Now, $g_{j_{m}}$ is in the span of our frame for each $m = 1,2,3,\ldots$,.  But, by inequality 3.1,
$$
\sum_{i=1}^{\infty}|<k_{i_{n}}g_{j_{m}}>|^{2} = \sum_{n=1}^{\infty}|k_{i_{n}}(j_{m})|^{2} < (\frac{1}{m})\|g_{j_{m}}\|.
$$
So this set is not a frame for its span.  This contradiction completes the proof of step I. 
  
\proclaim{Step II}
There is a natural number $m$ and numbers $A,B > 0$ so that $\forall j\in { \Omega}_{i}$, with $j\ge m$, we have:
$$
A  \le |h_{i}(j)| \le B.
$$ 
\endproclaim

To obtain the $m$, and the lower bound for $|h_{i}(j)|$, we proceed by way of contradiction.  If there is no such $m$ or $A$, then choose natural numbers $i_{1},j_{1}$ so that $0 < |h_{i_{1}}(j_{1})| \le 1$.  Since $h_{i_{1}}$ is finitely supported and for all $n\in \text{supp}h_{i_{1}}$, we have
$$
\sum_{i\in \Gamma}|h_{i}(n)|^{2} < \infty,
$$
it follows that there are natural numbers $i_{2} > i_{1}$, and $j_{2} > j_{1}$ with

(1)  $j_{2} > \text{max}\{\text{supp} h_{i_{1}}\}$ (so $h_{i_{1}}(j_{2}) = 0$),

(2)  $0 < |h_{i_{2}}(j_{2})| < \frac{1}{2}$,

(3)  $|h_{i_{2}}(j_{1})| \le \frac{1}{2}$.

Continuing by induction, we can find natural numbers $i_{1} < i_{2} < i_{3} < \cdots$ and $j_{1} < j_{2} < j_{3} < \cdots$ satisfying:

(4)  $h_{i_{n}}(j_{m}) = 0$, for all $m > n$,

(5)  $0 < |h_{i_{n}}(j_{n})| < \frac{1}{n}$,

(6)  $|h_{i_{n}}(j_{m})| \le \frac{1}{n}, \forall m < n$.

Let ${\Delta}^{c} = \{j_{i} : i = 1,2,3,\ldots\}$, and consider the subset of the frame $(f_{i})_{i=1}^{\infty}$ consisting of the elements, $((g_{i})_{i\in \Delta},(h_{i_{n}})_{n=1}^{\infty})$.  Let $P_{\Delta}$ be the natural orthogonal projection of $H$ onto $\text{span}(g_{i})_{i\in \Delta}$.  Note that (4)-(6) imply $\text{span}(g_{j_{k}})_{k=1}^{\infty} = \text{span}((I-P_{\Delta})h_{i_{n}})_{n=1}^{\infty}$.  By our assumption for this direction of the theorem, $((g_{i})_{i\in \Delta},(h_{i_{n}})_{n=1}^{\infty})$ is a frame sequence.  By Proposition 2.3 (2), $((I-P_{\Delta})h_{i_{n}})_{n=1}^{\infty}$ is also a frame sequence.   
Now, for all $n = 1,2,3,\ldots$, we have
$$
(I-P_{\Delta})h_{i_{n}} = \sum_{m=1}^{n}h_{i_{n}}(j_{m})g_{j_{m}}.
$$
But, 
$$
\text{inf}_{m}\sum_{n=1}^{\infty}|<(I-P_{\Delta})h_{i_{n}},g_{j_{m}}>|^{2} =
\text{inf}_{m}{\sum_{n=1}^{\infty}|h_{i_{n}}(j_{m})|^{2}} \le \text{inf}_{m}\sum_{n=m}^{\infty}|\frac{1}{n}|^{2} = 0, 
$$
which contradicts the fact that $((I-P_{\Delta})h_{i_{n}})_{n=1}^{\infty}$ is a frame for $\text{span}(g_{j_{m}})_{m=1}^{\infty}$.  This concludes the proof of step II.

Recall that $P_{m}$ denotes the natural (orthogonal) projection of $H$ onto $\text{span}(g_{i})_{i=1}^{m}$, and for $m < n$, $P_{m,n} = P_{n} - P_{m-1}$.

\proclaim{Step III}
There is a natural number $m_{0} > m$ so that $((g_{i})_{i=m_{0}+1}^{\infty},((I-P_{m_{0}})h_{i})_{i\in \Gamma})$ is a Riesz frame.
\endproclaim

We prove Step III by way of contradiction.  If no such $m_{0}$ exists, given m as in Step II, there are finite sets of natural numbers ${\Gamma}_{1}$  and ${\Delta}_{1} \subset \{n : n\ge m+1\}$ and a vector $f_{1} \in \text{span}\{(g_{i})_{i\in {\Delta}_{1}},((I - P_{m})h_{i})_{i\in {\Gamma}_{1}}\}$ satisfying:
$$
\|f_{1}\| = 1,
$$
$$
\sum_{i\in {\Delta}_{1}}|<f_{1},g_{i}>|^{2} + \sum_{i\in {\Gamma}_{1}}|<f_{1},(I - P_{m})h_{i}>|^{2} < 1.
$$
Let ${\ell}_{1} = \text{max}\{n: n\in {\Delta}_{1}\cup \cup_{i\in {\Gamma}_{1}}\text{supp}(I-P_{m})h_{i}\}$.  By Step II, there are only a finite number of $h_{i}$ whose supports intersect $\{m+1,m+2,\ldots,{\ell}_{1}\}.$  .  Since each $h_{i}$ has finite support, there is a natural number $m < m_{1}$ so that $(I - P_{m_{1}})h_{i} \not= 0$, implies $P_{m,{\ell}_{1}}h_{i} = 0$.  This fact and our assumption that Step III fails, implies the existence of finite sets of natural numbers ${\Gamma}_{2} \subset {{\Gamma}_{1}}^{c}$ and ${\Delta}_{2} \subset \{n: n\ge m_{1}\}$ and a vector $f_{2}$ satisfying: 
$$
h_{j} \in \text{span}\{(g_{n})_{n=m_{1}+1}^{\infty},(g_{n})_{n=1}^{m}\}, \ \ \ \ \forall j\in {\Gamma}_{2},
$$
$$
f_{2} \in \text{span}\{(g_{n})_{n\in {\Delta}_{2}},(I - P_{m_{1}+1})h_{i})_{i\in {\Gamma}_{2}}\}, \ \ \ \ \|f_{2}\| = 1,
$$ 

$$
\sum_{i\in {\Delta}_{2}}|<f_{2},g_{i}>|^{2} + \sum_{i\in {\Gamma}_{2}}|<f_{2},(I - P_{m_{1}})h_{i}>|^{2} < \frac{1}{2}.
$$

Continuing by induction, there exists natural numbers $m_{0} = m < m_{1} < m_{2} < \cdots$ and finite subsets of the natural numbers ${\Delta}_{i}$ and  ${\Gamma}_{i}$, and vectors $f_{i}$ satisfying:
$$
{\Delta}_{i} \subset \{m_{i-1}+1,m_{i-1}+2,\ldots,m_{i}\}, \tag 3.2
$$
$$
h_{j} \in \text{span}\{(g_{n})_{n=m_{i-1}+1}^{m_{i}},(g_{i})_{i=1}^{m}\}, \forall j\in {\Gamma}_{i}, \tag 3.3
$$
$$
f_{i}\in \text{span}\{(g_{n})_{n\in {\Delta}_{i}},(P_{m_{i-1},m_{i}}h_{j})_{j\in {\Gamma}_{i}}\}, \tag 3.4
$$
$$
\|f_{i}\| = 1,  \tag 3.5
$$
$$
\sum_{j\in {\Delta}_{i}}|<f_{i},g_{j}>|^{2} + \sum_{j \in {\Gamma}_{i}}|<f_{i},P_{m_{i-1},m_{i}}h_{j}>|^{2} < \frac{1}{i}.  \tag 3.6
$$

Next, let $\Delta = \cup_{i=1}^{\infty}{\Delta}_{i} \cup \{1,2,3,\ldots,m\}$ and $\Psi = \cup_{i=1}^{\infty}{\Gamma}_{i}$.  We will show that the subset of our frame given by $((g_{i})_{i\in \Delta},(h_{i})_{i\in \Psi})$ is not a frame for its closed linear span, contradicting our assumption that $(f_{i})_{i=1}^{\infty}$ has the subframe property.  To see this, let ${\Cal K} = \text{span}((g_{i})_{i\in \Delta},(h_{i})_{i\in \Psi})$ and note that $\{1,2,\ldots,m\} \subset \Delta$ and (3.3) imply that $P_{m_{i-1},m_{i}}h_{j} \in {\Cal K}$, for every $j\in {\Gamma}_{i}$.  Since ${\Delta}_{i} \subset \Delta$, it follows from (3.4) that $f_{i}\in {\Cal K}$, for all $i$.  Finally, by (3.2), (3.3), (3.4) we have:
$$
f_{i} \perp \text{span}((g_{j})_{j\in (\Delta - {\Delta}_{i})},(h_{j})_{j\in (\Psi - {\Gamma}_{i})}).
$$
Therefore,
$$
\sum_{j\in \Delta}|<f_{i},g_{j}>|^{2} + \sum_{j\in \Psi}|<f_{i},h_{j}>|^{2} =   
\sum_{j\in {\Delta}_{i}}|<f_{i},g_{j}>|^{2} + \sum_{j \in {\Gamma}_{i}}|<f_{i},h_{j}>|^{2} 
$$
$$
= \sum_{j\in {\Delta}_{i}}|<f_{i},g_{j}>|^{2} + \sum_{j \in {\Gamma}_{i}}|<f_{i},P_{m_{i-1},m_{i}}h_{j}>|^{2} < \frac{1}{i}.
$$
Therefore, our subset of the frame is not a frame sequence.  This completes the proof of Step III.

Now, let $G = \text{span}(g_{i})_{i=1}^{m_{0}}$, $P_{m_{0}}$ the natural (orthogonal) projection of $H$ onto $G$, and $h_{i}^{2} = P_{m_{0}}h_{i}$.  Also let  ${\Delta}_{i} = {\Omega}_{i}\cap \{m_{0}+1,m_{0}+2,\ldots\}$, and 
$$
h_{i}^{1} = (I - P_{m_{0}})h_{i} = \sum_{j\in {\Delta}_{i}}h_{i}(j), 
$$
 
\proclaim{Step IV}
We verify (ii), and (iii) of the theorem.
\endproclaim

Since $(h_{i})_{i\in \Gamma}$ is a frame and $\text{Rng}P_{m_{0}}$ is finite dimensional, (ii) follows from Lemma 3.1.  Part (iii) follows immediately from Step III and the fact that 
$$
\text{span}(g_{i})_{i=1}^{m_{0}} \perp \text{span}\{(g_{i})_{i=m_{0}+1}^{\infty},(h_{i}^{1})_{i\in \Gamma}\}
$$  

\vskip10pt

$(2) \Rightarrow (1)$  Assume that $((g_{i})_{i=1}^{\infty},(h_{i})_{i\in \Gamma},(k_{i})_{i=1}^{n})$ is a frame for $H$ satisfying the conditions in part (2) of the theorem.  Since $((g_{i})_{i=1}^{\infty},(h_{i}^{1})_{i\in \Gamma})$ is a Riesz frame, we assume it has the properties of Theorem 2.4.  Letting $\Delta$ equal the natural numbers and $D = \Gamma$ in (2)(ii) of Theorem 2.4, we get that $(h_{i}^{1})_{i\in \Gamma}$ is a frame sequence with lower frame bound $A_{0}$.  Since $\sum_{i\in \Gamma}\|h_{i}^{2}\|^{2} < \infty,$ there are only finitely many infinitely supported vectors $k_{i}$, and $(g_{i})_{i=1}^{\infty}$ is a Riesz basis it follows that our set of vectors satisfies the upper frame condition (and hence every subset satisfies the upper frame condition) with constant say $B$.  By taking the natural isomorphism of $(g_{i})_{i=1}^{\infty}$ to an orthonormal basis for $H$, we may assume that $(g_{i})_{i=1}^{\infty}$ is an orthonormal basis for $H$.  (To simplify the notation, we will use the same constants given earlier).  Choose an arbitrary subset of the frame of the form: $((g_{i})_{i\in \Delta},(k_{i})_{i\in \Lambda},(h_{i})_{i\in {\Gamma}_{1}})$.  Applying (2)(ii) of Theorem 2.4 again, we see that $(P_{{\Delta}^{c}}h_{i}^{1})_{i\in {\Gamma}_{1}}$ is a frame sequence with frame constants $A_{0},B$.  We will finish the proof in three steps.

\proclaim{Step I}
There is a subset $\Omega \subset {\Gamma}_{1}$ with $|{\Gamma}_{1}-\Omega|< \infty$, so that
$(P_{{\Delta}^{c}}h_{i})_{i\in \Omega}$ is a frame sequence.
\endproclaim

By our assumption (ii), we can choose $\Omega$ of the form above so that
$$
\sum_{i\in \Omega}\|h_{i}^{2}\|^{2} < (\frac{A_{0}}{2})^{2}.
$$
Then for any $f\in \text{span} 
(P_{{\Delta}^{c}}h_{i})_{i\in \Omega}$, and the fact that $(P_{{\Delta}^{c}}h_{i}^{1})_{i\in {\Gamma}_{1}}$ is a frame sequence with frame constants $A_{0},B$, we have
$$
\sqrt{\sum_{i\in \Omega}|<f,P_{{\Delta}^{c}}h_{i}>|^{2}} \ge \sqrt{\sum_{i\in \Omega}|<f,P_{{\Delta}^{c}}h_{i}^{1}>|^{2}} - \sqrt{\sum_{i\in \Omega}|<f,P_{{\Delta}^{c}}h_{i}^{2}>|^{2}} 
$$
$$
\ge \sqrt{A_{0}}\|f\| - \sqrt{\sum_{i\in \Omega}\|P_{{\Delta}^{c}}h_{i}^{2}\|^{2}\|f\|^{2}}
\ge \sqrt{A_{0}}\|f\| - \sqrt{\frac{A_{0}}{2}}\|f\| = \sqrt{\frac{A_{0}}{2}}\|f\|.  
$$ 

\proclaim{Step II}
We will prove that the family $((P_{{\Delta}^{c}}h_{i})_{i\in {\Gamma}_{1}},(P_{{\Delta}^{c}}k_{i})_{i\in \Lambda})$ is a frame sequence with frame bounds say $A,B$
\endproclaim

By step I, $(P_{{\Delta}^{c}}h_{i})_{i\in \Omega}$ is a frame sequence.  But $|\Lambda| < \infty$ and $|{\Gamma}_{1} - \Omega| < \infty$, and adding any finite number of vectors to a frame sequence always yields a frame sequence.

\proclaim{Step III}
We prove that $((g_{i})_{i\in \Delta},(h_{i})_{i\in \Gamma},(k_{i})_{i=1}^{n})$ is a frame sequence.
\endproclaim

Since $(g_{i})_{i=1}^{\infty}$ is an orthonormal basis, $P_{\Delta}$ is an orthoronal projection on $H$ with $I - P_{\Delta} = P_{{\Delta}^{c}}$.  Now, $(g_{i})_{i\in \Delta}$ is an orthonormal basis for its span, and by Step II we have that $((P_{{\Delta}^{c}}h_{i})_{i\in {\Gamma}_{1}},(P_{{\Delta}^{c}}k_{i})_{i\in \Lambda})$ is a frame sequence.  Applying Proposition 3.2 (1), it follows that $((g_{i})_{i\in \Delta},(h_{i})_{i\in \Gamma},(k_{i})_{i=1}^{n})$ is a frame sequence.  

This completes the proof of Theorem 3.2.          
\enddemo

Now, let us look at how this theorem uniquely relates frames with the subframe property to Riesz frames.  To get a frame with the subframe property, we first choose a Riesz frame $((g_{i})_{i=1}^{\infty},(h_{i})_{i\in \Gamma})$ for $H$ where $(g_{i})_{i=1}^{\infty}$ is a Riesz basis for $H$.  Now choose a finite set of vectors $(k_{i})_{i=1}^{n}$ from $H$ each with infinite support with respect to our Riesz basis $(g_{i})_{i=1}^{\infty}$.  Next, choose a natural number m and let $G = \text{span}(g_{i})_{i=1}^{m}$ be a finite dimensional subspace of $H$.    Finally, choose a set of vectors $(f_{i})_{i\in \Gamma}$ from $G$ satisfying:
$$
\sum_{i\in \Gamma}\|f_{i}\|^{2} < \infty.
$$
Then Theorem 3.2 yields that the set $((g_{i})_{i=1}^{\infty},(k_{i})_{i=1}^{n},(h_{i}+f_{i})_{i\in \Gamma})$ is a frame for $H$ which has the subframe property, and this is the only way to produce a frame with the subframe property.
This also shows, for example, that if we take a Riesz basis for $H$ and add to it an infinite number of infinitely supported vectors, then this new set has a subfamily which is not a frame for its closed linear span.

\heading{4.  The Projection Methods}
\endheading
\vskip10pt

If $(f_{i})_{i=1}^{\infty}$ is a frame, we define the {\bf frame operator} $S:H\rightarrow H$ by,
$$
S(f) = \sum_{i=1}^{\infty}<f,f_{i}>f_{i}. \tag 4.1
$$
Then $S$ is an isomorphism of $H$ onto $H$ and so $(S^{-1}f_{i})_{i=1}^{\infty}$ is also a frame for $H$ called the {\bf dual frame}.  For $f\in H$, we can write,
$$
f = SS^{-1}f = \sum_{i=1}^{\infty}<f,S^{-1}f_{i}>f_{i}, \tag 4.2
$$
where the $<f,S^{-1}f_{i}>$ are called the {\bf frame coefficients} for $f$.  One of the most difficult problems in frame theory is to explicitely calculate the dual frame of a frame.  A useful method here is to "truncate" the problem.  That is, for each $n$, let $H_{n} = \text{span}(f_{i})_{i=1}^{n}$ and $S_{n}:H_{n}\rightarrow H_{n}$ be given by,
$$
S_{n}f = \sum_{i=1}^{n}<f,f_{i}>f_{i}.  \tag 4.3
$$
For each $f\in H$, $P_{n}f$ converges to $f$ in norm.  But in general \cite{4}, the frame coefficients for $P_{n}f$ need not converge (even coordinatewise) to those of $f$.  If for every $f\in H$, and for every $i = 1,2,3,\ldots$, we have 
$$
\text{lim}_{n\rightarrow \infty}<f,S_{n}^{-1}f_{i}> = <f,S^{-1}f_{i}>, \tag 4.4
$$
we say that the {\bf projection method} works.  The advantage here is that finite dimensional methods, applied to the frame $(f_{i})_{i=1}^{n}$,  can be used to approximate the frame coefficients.  If $(<f,S_{n}^{-1}f_{i}>)_{i=1}^{n}$ converges to the frame coefficients for $f$ in the $\ell_{2}$ - sense, i.e.
$$
\text{lim}_{n\rightarrow \infty}\sum_{i=1}^{n}|<f,S_{n}^{-1}f_{i}> - <f,S^{-1}f_{i}>|^{2} + \sum_{i=n+1}^{\infty}|<f,S^{-1}f_{i}>|^{2} = 0, \tag 4.5
$$
we say that the {\bf strong projection method} works.  For a discussion of the projection method, we refer the reader to \cite{2}.  Also, for an in-depth study of the strong projection method, and a host of examples, we refer the reader to \cite{4}.  It is known \cite{2} that the projection method and the strong projection method working are not equivalent.  Also note that the projection methods depend upon the order in which the frame elements are written.  That is, a frame may satisfy the strong projection method but have a permutation which fails it \cite{4}.  It is immediate that the strong projection method works for Riesz bases (Or see Zwaan \cite{8}).  It also works for Riesz frames but may fail (even the projection method may fail) for frames with the subframe property \cite{4}.  The main theorem of this section will show that for frames with the subframe property, the projection methods become equivalent and independent of the order in which the frame elements are written.

\proclaim{Theorem 4.1}
If $(f_{i})_{i=1}^{\infty}$ is a frame with the subframe property, then the following are equivalent:

(1)  There are no infinitely supported vectors $k_{i}$ in Theorem 3.2,

(2)  $(f_{i})_{i=1}^{\infty}$ has a permutation satisfying the projection method,

(3)  Every permutation of $(f_{i})_{i=1}^{\infty}$ satisfies the strong projection method.
\endproclaim

\demo{Proof}
$(3) \Rightarrow (2)$ is obvious.

$(2) \Rightarrow (1)$  We will prove this by way of contradiction.  So suppose we have a frame $((g_{i})_{i=1}^{\infty},(h_{i})_{i\in \Gamma})$, $(k_{i})_{i=1}^{\ell})$ satisfying the conditions of Theorem 3.2.  As usual we may assume that $(g_{i})_{i=1}^{\infty}$ is an orthonormal basis for $H$.  Let $(f_{i})_{i=1}^{\infty}$ be a permutation of this frame satisfying the projection method.  Let $I,J$ be sets of natural numbers so that (recall the m of theorem 3.2):
$$
\{f_{i}:i\in I\} = \{g_{i}: 1\le i \le m\},  \ \ \ \ \{f_{i}:i\in J\} = \{k_{i}:1\le i \le \ell \}.
$$
Let $m_{0} = \text{max}_{i\in I\cup J}i$ and let $S_{n}$ be the frame operator for $(f_{i})_{i=1}^{n}$.  Our assumption that $(f_{i})_{i=1}^{\infty}$ satisfies the projection method implies there is a constant $K > 0$ so that for all $n\ge m_{0}$, we have $\|S_{n}^{-1}k_{i}\| \le K$.  So fix any $n\ge m_{0}$ and write $(f_{i})_{i=1}^{n}$ as $((g_{i})_{i\in \Delta}, (h_{i})_{i\in \Lambda}, (k_{i})_{i=1}^{\ell})$.  Let $Q_{n}$ be the orthogonal projection of $\text{span}(f_{i})_{i=1}^{n}$ onto its subspace $\text{span}\{(g_{i})_{i \in \Delta},(h_{i})_{i\in \Lambda}\}$.  Choose $1 \le j \le \ell$ so that
$$
\|(I - Q_{n})k_{j}\| = \text{max}_{1\le i \le \ell}\|k_{i}\|.  \tag 4.6
$$
Since the $h_{i}$ all have finite support with respect to the orthonormal basis $(g_{i})_{i=1}^{\infty}$, and the $k_{i}$ have infinite support, it follows that $\|(I - Q_{n})k_{j}\| \not= 0$ in formula (4.6).  Let 
$$
f_{n,j} = \frac{(I - Q_{n})k_{j}}{\|(I - Q_{n})k_{j}\|^{2}},
$$
so that $<f_{n,j},k_{j}> = 1$.  Finally, let
$$
f = f_{n,j} - \sum_{i\not= j}<f_{n,j},k_{i}>S_{n}^{-1}k_{i}.
$$
Now we compute,
$$
S_{n}f = S_{n}f_{n,j} - S_{n}(\sum_{i\not= j}<f_{n,j},k_{i}>S_{n}^{-1}k_{i}) 
$$
$$
= \sum_{i=1}^{n}<f_{n,j},f_{i}>f_{i} - \sum_{i\not= j}<f_{n,j},k_{i}>k_{i}
$$
$$
= \sum_{i=1}^{\ell}<f_{n,j},k_{i}>k_{i} - \sum_{i\not= j}<f_{n,j},k_{i}>k_{i}
= <f_{n,j},k_{j}>k_{j} = k_{j}.
$$
So $S_{n}^{-1}k_{j} = f$.  It follows from our earlier assumption that
$$
\|S_{n}^{-1}k_{j}\| = \|f\| \le K.  \tag 4.7
$$
Combining (4.6) with (4.7) we have
$$
K \ge \|f\| \ge \|f_{n,j}\| - \| \sum_{i\not= j}<f_{n,j},k_{i}>S_{n}^{-1}k_{i}\|  \tag 4.8
$$
$$
\ge \|f_{n,j}\| - \sum_{i\not= j}|<f_{n,j},(I - P)k_{i}>|\|S_{n}^{-1}k_{i}\|
$$
$$
\ge \|f_{n,j}\| - K\sum_{i\not= j}\|f_{n,j}\|\|(I - P)k_{i}\| \ge \|f_{n,j}\| - K{\ell}. 
$$
However, 
$$
\text{sup}_{n}\|f_{n,j}\| = \text{sup}_{n}\frac{1}{\|(I - Q_{n})k_{j}\|} = \infty.
\tag 4.9
$$
and (4.8) and (4.9) contradict one another.

$(1) \Rightarrow (3)$  By (1) and Theorem 3.2, our frame is of the form $((g_{i})_{i=1}^{\infty},(h_{i})_{i\in \Gamma})$ and has the properties listed in Theorem 3.2.  As usual, we may assume that $(g_{i})_{i=1}^{\infty}$ is an orthonormal basis for $H$.  Let $m$ and $G$ be given as in Theorem 3.2, and let $P_{G}$ be the natural (orthogonal) projection of $H$ onto $G$.  Let $(f_{i})_{i=1}^{\infty}$ be any permutation of this frame.  Choose a natural number $m_{0}$ so that $g_{j}\in \{f_{i}:1\le i \le m_{0}\}$, for all $1 \le j \le m$.  Let $A$ be the lower Riesz frame bound for the Riesz frame given in  Theorem 3.2 (iii), and choose $0 < \delta < \frac{1}{2}$ with 
$$
{\delta}^{2}\sum_{i\in \Gamma}\|P_{G}h_{i}\|^{2} < \frac{A}{4}.
$$
Let $n\ge m_{0}$ and let $S_{n}$ be the frame operator for $(f_{i})_{i=1}^{n}$.      
By our assumptions, there are finite sets of natural numbers $J\subset \Gamma$, and $I \subset \{m+1,m+2,\ldots$ so that $(f_{i})_{i=1}^{n} = ((g_{i})_{i=1}^{m},(g_{i})_{i\in I},(h_{i})_{i\in J})$.  Choose $f\in  \text{span}(f_{i})_{i=1}^{n}$ with 
$$
1 = \|f\|^{2} = \|P_{G}f\|^{2} + \|(I - P_{G})f\|^{2}.
$$
We consider two cases.

\proclaim{Case I}
$\|P_{G}f\|^{2} \ge \delta$.
\endproclaim

In this case,
$$
\sum_{i=1}^{n}|<f,f_{i}>|^{2} \ge \sum_{i=1}^{m}|<f,g_{i}>|^{2} = \|P_{G}f\|^{2} \ge \delta.
$$

\proclaim{Case II}
$\|P_{G}f\|^{2} \le \delta$.
\endproclaim

In this case, applying (iii) of Theorem 3.2, we have,
$$
\sqrt{\sum_{i=1}^{n}|<f,f_{i}>|^{2}} \ge \sqrt{\sum_{i\in I}|<f,g_{i}>|^{2} + \sum_{i\in J}|<f,h_{i}>|^{2}}   
$$
$$
\ge \sqrt{\sum_{i\in I}|<(I-P_{G})f,g_{i}>|^{2} + \sum_{i\in J}|<(I-P_{G})f,(I-P_{G})h_{i}>|^{2}} -
$$
$$
\sqrt{\sum_{i\in J}|<P_{G}f,P_{G}h_{i}>|^{2}}
$$
$$
\ge \sqrt{A}\|(I-P_{G})f\| - \sqrt{\sum_{i\in J}\|P_{G}f\|^{2}\|P_{G}h_{i}\|^{2}}
$$
$$
\ge \sqrt{{A}(1-\delta)} - {\delta}\sqrt{\sum_{i\in J}\|P_{G}h_{i}\|^{2}}
\ge \sqrt{\frac{A}{2}} - \sqrt{\frac{A}{4}}.
$$
Hence, our frame $(f_{i})_{i=1}^{\infty}$ satisfies the strong projection method.     
\enddemo

Although frames with the subframe property may fail even the projection method, Theorem 4.1 implies that this occurs because of a few "misbehaved" vectors.  We state this formally as,

\proclaim{Corollary 4.2}
If $(f_{i})_{i\in I}$ is a frame with the subframe property, then there is a finite subset $\Delta \subset I$ so that the strong projection method works for $(f_{i})_{i\in I-\Delta}$. 
\endproclaim

\proclaim{ACKNOWLEDGEMENT}
The author thanks Ole Christensen for helpful discussions concerning the material in this paper.
\endproclaim


\Refs

\ref\no{1}
\by  Berberian, Sterling K.  
\paper Introduction to Hilbert Space  
\jour University Texts, Cambridge University Press
\yr 1961
\endref

\ref\no{2}
\by  Casazza, Peter G. and Ole Christensen  
\paper  Riesz frames and Approximation of the Frame Coefficients
\jour (preprint)
\endref



\ref\no{3}
\by  Casazza, Peter G. and Ole Christensen  
\paper Hilbert Space Frames Containing a Riesz basis and Banach Spaces Which Have No Subspace Isomorphic To $c_{0}$.   
\jour  (to appear, Journal of Math. Analysis and Applications)
\endref

\ref\no{4}
\by  Christensen, Ole 
\paper Frames and the Projection Method  
\jour Applied and Computational Harmonic Analysis
\vol 1
\yr 1993
\pages 50-53
\endref




\ref\no{5}
\by  Christensen, Ole  
\paper Frames and Pseudo-Inverses  
\jour Journal of Math. Analysis and Applications
\vol 195
\yr 1995
\pages 401-414
\endref




\ref\no{6}
\by  Christensen, Ole  
\paper Frames Containing a Riesz Basis and Approximation of the Frame Coefficients using Finite-Dimensional Methods  
\jour  J. Math. Anal. and Appl.
\vol 199
\yr  1996
\pages  256-270
\endref

\ref\no{7}
\by Young, R.M.   
\paper An Introduction to Nonharmonic Fourier Series  
\jour Academic Press, New York
\yr 1980
\endref


\ref\no{8}
\by  Zwaan, M.  
\paper Approximation of the Solution to the Moment Problem in a Hilbert Space  
\jour Numerical Functional analysis and Optimization
\vol 11
\yr 1990
\pages 601-608
\endref

\endRefs
\enddocument